\documentclass[12pt]{article}

\usepackage[a4paper]{geometry}
\usepackage[myheadings]{fullpage}
\usepackage{fancyhdr}
\usepackage{lastpage}
\usepackage{graphicx, wrapfig, subcaption, setspace, booktabs}
\usepackage[T1]{fontenc}
\usepackage[font=small, labelfont=bf]{caption}
\usepackage[protrusion=true, expansion=true]{microtype}
\usepackage[english]{babel}
\usepackage{sectsty}
\usepackage{url, lipsum}
\usepackage{tgbonum}
\usepackage{xcolor}
\usepackage{tikz}
\usepackage{cancel}
\usepackage{lineno,hyperref}

\usepackage[utf8]{inputenc}
\usepackage{amsmath,amssymb,amsthm }
\usepackage{framed}
\usepackage{varioref, cleveref}
\usepackage{enumerate}
\usepackage{lmodern}
\usepackage{tikz}
\usepackage{pdfpages}
\usepackage{float}
\usepackage{pdfcomment}
\usepackage{verbatim}
\usepackage{multirow}

\usepackage[ruled,vlined]{algorithm2e} %
\newcommand{\F}{\mathcal{F}}
\newcommand{\Dtau}{\Delta \tau}

\newcommand{\tY}{\Tilde{Y}}
\newcommand{\half}{\frac{1}{2}}

\newcommand*{\bm}[2]{W_{#1}^{#2}}  

\newcommand{\pa}{\partial}

\newcommand{\pd}[2]{\frac{\pa {#1}}{\pa {#2}}}

\newcommand{\Vd}[0]{V_\text{data}}
\newcommand{\xir}[0]{\xi_\text{ref}}
\newcommand{\sn}[0]{\sigma_\nu}
\newcommand{\kn}[0]{\kappa_\nu}
\newcommand{\mr}[0]{\mu_{\nu}}
\newcommand{\tV}[0]{\tilde{V}}
\newcommand{\snt}[0]{\tilde{\sigma}_\nu}
\newcommand{\knt}[0]{\tilde{\kappa}_\nu}
\newcommand{\mrt}[0]{\tilde{\mu}_{\nu}}
\newcommand{\rhot}[0]{\tilde{\rho}}
\newcommand{\tb}[0]{\tilde{b}}
\newcommand{\tOmega}[0]{\tilde{\Omega}}

\newcommand{\dtau}[1]{\pd{{#1}}{\tau}}

\floatname{algorithmic}{Procedure}

\begin{document}
\title{A gradient based calibration method \\ for the Heston model}
\author{Anna Clevenhaus, Claudia Totzeck and Matthias Ehrhardt}

\maketitle
\begin{abstract}
The Heston model is a well-known two-dimensional financial model. 
Because the Heston model contains implicit parameters that cannot be determined directly from real market data, calibrating the parameters to real market data is challenging. 
In addition, some of the parameters in the model are non-linear, 
which makes it difficult to find the global minimum of the optimization problem within the calibration. 
In this paper, we present a first step towards a novel space mapping approach for parameter calibration of the Heston model. 
Since the space mapping approach requires an optimization algorithm, we focus on deriving a gradient descent algorithm. 
To this end, we determine the formal adjoint of the Heston PDE, which is then used to update the Heston parameters.  
Since the methods are similar, we consider a variation of constant and time-dependent parameter sets.
Numerical results show that our calibration of the Heston PDE works well for the various challenges in the calibration process and meets the requirements for later incorporation into the space mapping approach. 
Since the model and the algorithm are well known, this work is formulated as a proof of concept.
\end{abstract}

\section{Introduction}
In finance, calibrating model parameters to fit real market data is challenging because most model parameters are implicit in the real market data \cite{Cui.2017,horvath2021deep,liu2019neural,mehrdoust2023two,Mikhailov.2003,Teng.2016}. 
We consider the well-known two-dimensional Heston model, which contains at least four parameters implicit in the market data. 
The Heston proposed a two-dimensional stochastic differential equations (SDE) model to simulate the behavior of the stock price \cite{Heston.1993} and presented a closed-form valuation formula for his model. 
Some calibration techniques are based on this formula \cite{Mikhailov.2003,Cui.2017}. 
The closed-form equation has some restrictions w.r.t.\ improving the model by considering non-constant parameters \cite{Mikhailov.2003}.
Another strategy to improve the model is to introduce additional processes, which leads to an increase in the difficulty of the calibration process \cite{Teng.2016}. 
Since the stock price and variance are stochastic processes in the Heston model, one can  
Monte Carlo optimization methods can be applied \cite{Teng.2016}.

We consider the well-known two-dimensional Heston model, which contains at least four parameters implicit in the market data. 
The parameter calibration is formulated as a constrained optimization problem to minimize a cost functional. 
The cost functional describes the difference between the reference data, the subsequent market data, and the data obtained by numerically solving our model. 
In the space mapping approach \cite{bandler2004space}, a coarse and a fine optimization solver are used.
This paper presents the first step "a gradient descent algorithm" for the Heston model towards a space mapping approach in finance. 
Later, within space mapping, the gradient descent algorithm will compute a coarse and a cheap approximation for the calibration problem of the Heston model. 

Our goal is to introduce the space mapping approach to financial applications, and therefore we use the $\log$-transformed normalized Heston model formulated as a partial differential equation (PDE) and the gradient descent algorithm as a pre-step, since both are well known and thus we can focus on the novel aspects. 
During our research, we didn't find any other PDE-based calibration approach for the Heston model.
Of course, we are aware that there are already algorithms to compute the exact solution of the Heston calibration problem \cite{Cui.2017}. 
However, these methods are limited to the assumption that the parameters are constant, whereas our approach considers time-dependent parameters and constant parameter calibration. 
Furthermore, if the parameters are assumed to be time-dependent, an analytical solution can't be derived and thus an exact solution is not available.
Therefore, as mentioned before, the paper is considered as a proof of concept to introduce a new calibration method based on financial research. 

We focus on deriving a gradient descent algorithm for the Heston model that satisfies the requirements for its use in the space mapping approach. 
In addition to our financial model, the gradient descent algorithm requires the adjoint of the model. 
Therefore, we formally derive the adjoint of the Heston model and construct the gradient using well-known techniques from optimization of partial differential equations (PDEs) \cite{HUUP,Troeltzsch}. 
The gradient descent algorithm has previously been applied to the Heston model using neural networks, so it has not been explicitly computed \cite{Liu.2019}. 
In this work, we focus on the calibration rather than the numerical solution of the model, i.e., our obtained results can be further improved by using more accurate numerical methods. Nevertheless, the calibration results obtained are remarkable.
Our attempt can be easily extended to include time-dependent parameters for the variance, thus providing an improvement over current calibration strategies based on the semi-analytical solution of the Heston model.

The rest of the article is organized as follows.
In Section~\ref{sec: Financial modelling}, we introduce the Heston model and our approach to solving the resulting $\log$-transformed Heston PDE. 
We then focus on parameter calibration. 
In Section~\ref{sec: Gradient Descent Algorithm}, we derive the adjoint of the $\log$-transformed Heston model and the corresponding gradient. 
We then discretize the method in Section~\ref{sec: Discretization} 
Finally, in Section~\ref{sec: Numerical Results} we focus on numerical results for the application of the gradient descent algorithm.
We analyze the behavior of the algorithm in various problems for constant and time-dependent parameter calibration. 
The paper ends with a conclusion and an outlook in Section~\ref{sec: Conclusion and Outlook}.

\section{The Heston Model}\label{sec: Financial modelling}
The Heston model was proposed by Heston in 1993 \cite{Heston.1993} and describes the dynamics of the underlying asset $S\in \mathbb{R}^{>0}$ and the variance $\nu\in \mathbb{R}^{>0}$ by a two-dimensional SDE. 
This model is an extension of the well-known Black-Scholes model, %
which considers only a stochastic process for the asset and leads to
\begin{equation}
    dS_t = r S_t\,dt + \sigma_S S_t\,d\bm{t}{S},\quad S_0>0,\\
\end{equation}
where $r\in \mathbb{R}^{>0}$ is the risk-free rate, $\sigma_S\in \mathbb{R}^{>0}$ is the volatility, and $\bm{t}{S}$ is the Brownian motion. 
By definition, the variance is the square of the volatility of the asset, $\nu=\sigma_S^2$. 
Heston considered a Cox-Ingersoll process for modeling the variance, which leads to the SDE system of Heston's model under the risk-neutral measure given by
\begin{equation}\label{eq:pureHeston}
 \begin{cases}
 dS_t = r S_t\,dt + \sqrt{\nu_t} S_t\,d\bm{t}{S},&\quad S_0>0,\\
 d\nu_t = \kn(\mr-\nu_t)\,dt+\sn \sqrt{\nu_t}\,d\bm{t}{\nu},&\quad \nu_0>0,
\end{cases},
\end{equation}
where $\kn\in\mathbb{R}^{>0}$ is the mean reversion rate, $\mr\in \mathbb{R}^{>0}$ is the long term mean, and 
and $\sn\in\mathbb{R}^{>0}$ is the volatility of variance. 
The Brownian motions $\bm{t}{S}$ and $\bm{t}{\nu}$ are correlated by the constant $\rho \in [-1,1]$ \cite{Heston.1993}.
For the variance process to be positive, the Feller condition $2\kn \mr \ge\sn^2$ must be satisfied. 
If the Feller condition is violated, the \eqref{eq:pureHeston} becomes complex which leads to computational problems.
Using Kolmogorov's backward equation, we derive the Heston PDE under the risk-neutral measure
\begin{equation}\label{eq: Heston transformed_notlog}
\begin{split}
    0= V_t + \frac{\nu}{2} S^2 V_{SS}+ \frac{1}{2} \sigma_{\nu}^2 \nu V_{\nu \nu}+ rS V_S+\kn (\mr-\nu)V_{\nu}+ \sigma_{\nu} \nu \rho S V_{S \nu}-rV ,\\
\end{split}
\end{equation}
where $V(S,\nu,t)$ denotes the fair price of the option.
If we look at the European plain vanilla put option, the terminal condition ("pay-off") is as follows
\begin{equation}\label{eq:payoff}
    V(S,\nu,T)=(K-S,0)^+,
\end{equation}
where $K\in \mathbb{R}^{>0}$ denotes the predefined strike price.
Heston proposed the following boundary conditions 
\begin{align}
    V(0,\nu,t)&=K\exp(-r(T-t)),\label{eq: Boundary Heston 1}\\
    V(S,\nu,t)&=0,\quad\text{for }S\to\infty,\label{eq: Boundary Heston 2}\\
    V(S,0,t)&=rs V_S(S,0,t)+ \kn \mr V_{\nu}(S,0,t)-rV(S,0,t),\label{eq: Boundary Heston 3} \\
    V(S,v,t)&\sim K\exp(-r(T-t)),\quad\text{for }\nu\to\infty.\label{eq: Boundary Heston 4}
\end{align}
The boundary conditions for the underlying follow directly from the assumptions on the financial market. 
We reverse the time direction by introducing $\tau=T-t$ and thus the payoff condition \eqref{eq:payoff} becomes an initial condition. 
Next, it is advantageous to use the variable transformation $x=\log(S/K)$ for the asset, as it simplifies the numerical implementation, i.e., we obtain a $\log$-transformed variable with normalization to the strike price $K$. 
We rewrite the option price equation as
$\tV(x,\nu,\tau)=V(S,\nu,t)/K$ and obtain the so called \textit{$\log$-transformed normalized Heston PDE}
\begin{equation}\label{eq: Heston transformed}
\begin{split}
    \tV_{\tau}= \frac{\nu}{2} \tV_{xx}+ \frac{1}{2} \sigma_{\nu}^2 \nu \tV_{\nu \nu}+ (r-q-\frac{\nu}{2})\tV_x+\kn (\mr-\nu)\tV_{\nu}+ \sigma_{\nu} \nu \rho \tV_{x \nu} -r\tV ,
\end{split}
\end{equation}
defined on the semi-unbounded domain $x\in\mathbb{R}$, $\nu\ge0$, $0\le\tau\le T$ and supplied with the initial condition
\begin{equation}
    \tV(x,\nu,\tau)=\bigl(1-\exp(x),0\bigr)^+
\end{equation}
and the following boundary conditions
\begin{align}
    \tV(x,\nu,\tau)&\sim\exp(-r\tau), \quad\text{for }x\to-\infty,\\
    \tV(x,\nu,\tau)&=0,\quad\text{for }x\to\infty,\\
    \tV(x,\nu,\tau)&\sim\exp(-r\tau),\quad\text{for }\nu\to\infty.\label{bc_nu_infty}
\end{align}
The boundary conditions for the underlying asset are the normalized $\log$-transformation of Heston's proposed boundary conditions \eqref{eq: Boundary Heston 1}, \eqref{eq: Boundary Heston 2} and \eqref{eq: Boundary Heston 4}.
At $\nu=0$ the parabolic PDE degenerates to a first order hyperbolic PDE, and therefore we need to consider the Fichera theory \cite{Buckova.2016,Kutik.2015} to assess whether it is necessary to provide these analytic boundary conditions or not.  

From the Fichera condition at $\nu=0$ given by 
\begin{equation}
   b(\nu)= \kn ( \mr-\nu) - \half \sn^2 
\end{equation}
it follows
\begin{itemize}
    \item if $\lim_{\nu \to 0^+} b(\nu)\ge0$ (outflow) we must not supply any boundary conditions at $\nu=0$.
    \item if $\lim_{\nu \to 0^+} b(\nu)<0$ (inflow) we have to supply boundary conditions at $\nu=0$.
\end{itemize}
In the sequel we assume that the Feller condition $2\kn\mr\ge\sn^2$ is satisfied. 
Therefore we obtain an outflow boundary at $\nu=0$ and must not impose any analytical boundary conditions on this boundary.  
But we obviously need a numerical boundary condition to complete the scheme, which will be later discussed in section \ref{sec: Discretization}.

In addition to the Heston PDE with constant parameters for the variance process, we also consider the time-dependent parameters $\knt$, $\mrt$, $\snt$, as well as $\rhot$. Then the corresponding PDE formulation is given by
\begin{equation}\label{eq: Heston transformed time}
    \tV_{\tau}= \frac{\nu}{2} \tV_{xx}+ \frac{1}{2} \snt^2 \nu \tV_{\nu \nu}+ (r-q-\frac{\nu}{2})\tV_x+\knt (\mrt-\nu)\tV_{\nu}+ \snt \nu \rhot \tV_{x \nu}-r\tV.
\end{equation}
In the numerical simulations \ref{sec: Numerical Results}, we will also discuss variations of constant and time-dependent parameter.

\section{Derivation of a gradient based optimization strategy}\label{sec: Gradient Descent Algorithm}
For a given data set $\Vd$ and reference parameters $\xir$, we formulate the calibration as optimization task with the cost functional given by
\begin{equation*}
   J(\tV,\xi) = \frac{1}{2} \int_0^T \| K \tV - \Vd \|^2 \,d\tau + \frac{\lambda}{2} \| \xi - \xir \|^2.
\end{equation*}
As our aim is to fit the parameter to real market data and therefore no $\xir$ exists, we set $\lambda=0$ and the cost functional reduces to 
\begin{equation}\label{eq: Cost functional}
   J(\tV,u) = \frac{1}{2} \int_0^T \| K \tV - \Vd \|^2 \,d\tau.
\end{equation}
In the following, we derive a gradient-based algorithm that allows us to calibrate the parameters numerically. To this end we use a Lagrangian approach.

\subsection{First-order optimality conditions for the Heston model}
Let us denote the Lagrange multipliers by $\psi = (\varphi, \varphi^a, \varphi^b, \varphi^c,\varphi^d)$, set $\Omega = (0,\nu_{\max}) \times (x_{\min}, x_{\max})$ and split the boundary $\partial \Omega$ into 
\begin{align*}
 &\Gamma_a = \partial\Omega \cap \{x=-\infty\}, && \Gamma_b = \partial\Omega \cap \{x=\infty\}, \\
 &\Gamma_c = \partial\Omega \cap \{\nu=0\}, && \Gamma_d = \partial\Omega \cap \{\nu=\infty\}.
\end{align*}
First, we focus on the $\log$-transformed normalized Heston equation with constant parameters \eqref{eq: Heston transformed}  and define
\begin{align*}
A = \frac{1}{2} \nu \begin{pmatrix} \sn^2 & \sn \rho \\ \sn \rho & 1  \end{pmatrix}, \qquad 
b=\begin{pmatrix} \kn(\mr - \nu) -  \frac{1}{2} \sn^2 \\r-q-\frac{\nu}{2} - \frac{1}{2} \sn \rho \end{pmatrix} \qquad 
\tb=\begin{pmatrix} \kn\mr \\r-q \end{pmatrix}.
\end{align*}
Then it can be written as
\begin{equation*}
    \tV_\tau - \nabla \cdot A \nabla \tV - b \cdot \nabla \tV + r\tV = 0.
\end{equation*}
Next, we define the operator $e$ which will represent the constraint in the Lagrangian. 
Since at $\Gamma_c$ no boundary condition has to be given, we introduce $\tilde{\Omega}=\Omega \cap \Gamma$ and the operator $e$ is implicitly defined by 
\begin{equation} \label{eq: Lagrangian constraint}
\begin{aligned}
    \langle e(\tV,\xi) , \psi \rangle &= \int_0^T \int_{\tOmega} \Bigl[ \tV_{\tau} - \nabla \cdot A \nabla \tV - b \cdot \nabla \tV + r\tV \Bigr] \varphi\, dz d\tau \\ 
    &\qquad + \int_0^T \int_{\Gamma_a} \bigl[ \tV - \exp(-r\tau) \bigr] \varphi^a \,ds d\tau + \int_0^T \int_{\Gamma_b} \tV  \varphi^b\, ds d\tau \\
    &\qquad + \int_0^T \int_{\Gamma_d} \bigl[  \tV - \exp(-r\tau)\bigr] \varphi^d\, ds d\tau \\
      &=: T_1 + T_2 + T_3 + T_4. 
\end{aligned}
\end{equation} 

The Lagrangian for the constrained parameter calibration problem is then given by
\begin{equation*}
   L(\tV,\xi,\psi) = J(\tV,\xi) - \langle e(\tV,\xi), \psi \rangle.
\end{equation*}
We formally compute the first-order optimality conditions by setting $dL = 0$.
For details on the method we refer to \cite{HUUP,Troeltzsch}.
Before computing the G\^ateaux derivatives of $L$ in arbitrary directions \cite{HUUP}, we note that by Green's first identity it holds  
\begin{equation}\label{eq:HestonX}
   \int_\Omega (b\cdot \nabla \tV) \varphi \,dz = \int_{\partial \Omega} (b \cdot \vec n) \tV \varphi \,ds - \int_\Omega \tV \nabla \cdot (b\varphi) \,dz.
\end{equation}
Therefore, we can rewrite 
\begin{align*}
T_1 &= \int_0^T \int_{\tOmega} \varphi \tV_\tau + A \nabla \tV \cdot \nabla \varphi - \frac{1}{2} b\cdot \nabla \tV \varphi +\frac{1}{2} \tV b \cdot \nabla \varphi + (r + \frac{1}{2} \nabla \cdot b) \tV \varphi \,dz \\
&\qquad - \int_{\partial \Omega} (A \nabla \tV ) \cdot \vec n \varphi \,ds 
   - \frac{1}{2} \int_{\partial\Omega} (b\cdot \vec n) \tV \varphi\, ds d\tau\\
&=\left[ \int_\Omega \varphi \tV \,dz\right]_{\tau=0}^{\tau=T} 
    + \int_0^T \int_{\tOmega} \tV \Bigl[- \varphi_\tau - \nabla \cdot A^\top \nabla \varphi + b\cdot \nabla \varphi + (r+ \nabla\cdot b) \varphi \Bigr] \,dz   \\
&\qquad+ \int_{\partial \Omega} [ (A^\top\nabla \varphi) \cdot \vec n  - (b\cdot \vec n) \varphi]\tV \,ds
      - \int_{\partial \Omega} (A \nabla \tV ) \cdot \vec n \varphi \,ds d\tau.
\end{align*}
As $e$ is linear in $V$, the G\^ateaux derivative in some arbitrary direction $h$ reads
\begin{align*}
d_V \langle e(V,\xi),\psi \rangle [h] 
&= 
  \left[ \int_\Omega \varphi h \,dz\right]_{\tau=0}^{\tau=T} 
  + \int_0^T \int_\Omega h \bigl[ -\varphi_\tau- \nabla \cdot A^\top \nabla \varphi + b\cdot \nabla \varphi + (r+ \nabla\cdot b) \varphi \bigr] dz   \\
&\qquad+ \int_{\partial\Omega} \bigl[ (A^\top\nabla \varphi) \cdot \vec n  - (b\cdot \vec n) \varphi \bigr]h \,ds 
   - \int_{\partial \Omega} (A \nabla h ) \cdot \vec n \varphi \,ds\, d\tau \\
&\qquad + \int_0^T \int_{\Gamma_a} h \varphi^a \,ds\, d\tau 
       + \int_0^T \int_{\Gamma_b} h  \varphi^b\, ds\, d\tau 
       + \int_0^T \int_{\Gamma_d} h  \varphi^d\, ds\, d\tau.
\end{align*}
For the cost functional we have
\begin{equation*} 
    d_{\tV} J(\tV,\xi) [h] = \int_0^T \int_\Omega h K (K\tV-\Vd) \,dz \,d\tau.
\end{equation*}

To identify the adjoint equation, we consider
\begin{equation*}
    \begin{aligned}
    0 &=d_{\tV} L(\tV,u,\psi) [h] \\
    &= \int_0^T h \left[ \int_\Omega K(K\tV-\Vd) + \varphi_\tau + \nabla \cdot A^\top \nabla \varphi - b\cdot \nabla \varphi - (r+ \nabla\cdot b) \varphi \right] \, dz\, d\tau.
    \end{aligned}
\end{equation*}
for arbitrary $h$.
Note that we are not allowed to vary $\tV$ at $\tV(x,\nu,0)$ as the initial condition is fixed. Therefore we have $h(0,\nu,x)\equiv0$.

For choosing $h\equiv0$ on $\partial \Omega$ and $h(T,\nu,x)=0$, we find with the Variational Lemma 
\begin{equation*}
  \varphi_\tau + \nabla \cdot A^\top \nabla \varphi - b\cdot \nabla \varphi - (r+ \nabla\cdot b) \varphi = -K(K\tV-\Vd) \quad \text{on }\Omega.
\end{equation*}
Now, choosing $h(T,\nu,x) \ne 0$, we then obtain the terminal condition $\varphi(T,\nu,x) = 0$.

We consider the four boundary conditions separately. 
At $\Gamma_c$, also the parabolic adjoint PDE degenerates to a first-order hyperbolic PDE, 
and thus we have to consider the Fichera theory \cite{Buckova.2016,Kutik.2015} for the variance again.

The Fichera condition w.r.t.\ the variance at $\nu=0$ of the adjoint is the same as before. 
Therefore no analytic boundary condition is supplied for this boundary, as we assume that the Feller condition holds.
On $\Gamma_a$ we have
\begin{equation}\label{eq:GammaA}
 0=\int_{\Gamma_a} \bigl[ (A^\top\nabla \varphi) \cdot \vec n  - (b\cdot \vec n) \varphi\bigr]h  - (A \nabla h ) \cdot \vec n \varphi  + h \varphi^a \,ds.
\end{equation}
Choosing $h\equiv\text{const}\ne 0$ yields
 \begin{align*}
 0&= \int_{\Gamma_a}  h \bigl[ (A^\top\nabla \varphi) \cdot \vec n  - (b\cdot \vec n) \varphi  +  \varphi^a \bigr] \,ds,
\end{align*}
hence $(A^\top\nabla \varphi) \cdot \vec n  - (b\cdot \vec n) \varphi  +  \varphi^a =0$.
On the other hand, choosing $\nabla h \ne 0$ \eqref{eq:GammaA} must still hold. 
This yields $\varphi = 0$ on $\Gamma_a$ and 
\begin{align*}
\varphi^a=-(A^\top\nabla \varphi) \cdot \vec n
=-(A^\top\nabla \varphi) \cdot \begin{pmatrix} 0 \\-1 \end{pmatrix}
=\half \nu \sn\rho \varphi_{\nu} + \half \nu \varphi_x=\half \nu \varphi_x
\end{align*}
As $\tV(x_{\min},\nu,0)=\exp(-r\tau)$ is given and independent of $x$, we obtain $\varphi_x=0$ there and thus $\varphi^a=\varphi=0$ at this boundary.
Similarly, we find on $\Gamma_{b}$ that
\begin{equation}\label{eq:GammaB}
    0=\int_{\Gamma_b} \bigl[ (A^\top\nabla \varphi) \cdot \vec n  - (b\cdot \vec n) \varphi\bigr]h  - (A \nabla h ) \cdot \vec n \varphi  + h \varphi_{b}\,ds.
\end{equation}
With the same arguments, we obtain $\varphi = 0$ and $\varphi^b=0$ on $\Gamma_b$ as well.
Following the same arguments for $\Gamma_d$ we obtain $(A^\top\nabla \varphi) \cdot \vec n  - (b\cdot \vec n) \varphi  +  \varphi^d =0$ with $\varphi = 0$ and thus it reduces to
\begin{align*}
\varphi^d=-(A^\top\nabla \varphi) \cdot \vec n
=-(A^\top\nabla \varphi) \cdot \begin{pmatrix} 1 \\0 \end{pmatrix}
=-\half \nu \sn^2\varphi_{\nu} + \half \nu \sn \rho \varphi_x=-\half \nu \sn^2\varphi_{\nu}.
\end{align*}
As $\tV(x,\infty,\tau)=\exp(-r\tau)$ is given and independent of $\nu$, we obtain $\varphi_{\nu}=0$ there and thus $\varphi^d=\varphi=0$ at this boundary.

Altogether, the adjoint equation reads
\begin{equation}\label{eq: Heston AD}
    \dtau \varphi + \nabla \cdot A^\top \nabla \varphi - b\cdot \nabla \varphi 
    - (r+ \nabla\cdot b) \varphi = -K(K\tV-\Vd) \quad \text{on }\Omega,
\end{equation}
which is equivalent to
\begin{multline}\label{eq:HestonAD_noA}
     \varphi_{\tau} + \frac{1}{2} \nu \sn^2 \varphi_{\nu \nu}+ \nu \sn\rho \varphi_{x \nu} + \frac{1}{2} \nu \varphi_{xx}+ (\sigma_\nu^2-\kappa_\nu(\mr - \nu) ) \varphi_{\nu} \\
     + (q-r+\frac{\nu}{2} +\sigma_\nu \rho) \varphi_{x} + (\kn -r) \varphi = -K(K\tV-\Vd) \quad \text{on }\Omega
\end{multline}
with terminal condition $\varphi(T)=0$ and $\varphi=0$ on the boundaries $\Gamma_a$, $\Gamma_b$ and $\Gamma_d$ and the outflow boundary at $\nu=0$.

\subsection{Derivation of the Gradient}
Let $\xi=(\sn,\rho,\kn,\mr)$ be the parameters to be identified, as $r$ and $q$ are given by the data.
We compute the optimality condition by setting $d_{\xi} L(V,\xi,\psi) = 0$. 
Since the boundaries $\Gamma_a$, $\Gamma_b$ and $\Gamma_d$ are zero, we focus on $\tOmega$. 
In the following we state the derivatives w.r.t.\ the different parameters separately.
For $\sn$ we get
\begin{equation*}
    d_{\sn} \langle e(\tV,\xi),\psi \rangle  = \int_0^T \int_{\tilde{\Omega}}  \tV \Bigl[ - \sn \nu \varphi_{\nu \nu} - 2 \sn \varphi_{\nu} -\rho \varphi_x - \rho \nu \varphi_{x \nu}\Bigr]\, dz\, d\tau.
    \end{equation*}
Similarly, we obtain for the other derivatives
\begin{align*}
    d_{\rho} \langle e(\tV,\xi),\psi \rangle  
    &=   \int_0^T \int_{\tilde{\Omega}} \tV \Bigl[ - \sn \varphi_x - \sn \nu \varphi_{x \nu}\Bigr]\,dz\, d\tau, \\
    d_{\kn} \langle e(\tV,\xi),\psi \rangle 
    &= \int_0^T \int_{\tilde{\Omega}} \tV \Bigl[ (\mr-\nu)\varphi_{\nu}-\varphi \Bigr] \,dz\, d\tau, \\
     d_{\mr} \langle e(\tV,\xi),\psi \rangle 
     &=  \int_0^T \int_{\tilde{\Omega}}  \kn \tV \varphi_{\nu} \, dz\, d\tau.
\end{align*}
Note that $d_{\xi} L(\tV,\xi,\psi)[h_{\xi}] = 0$ needs to hold for arbitrary directions $h_{\xi}$. 
Therefore, we can read off the gradient from the above expressions.

We extend this gradient formulation for time dependent parameter $\tilde{\xi}=(\snt,\rhot,\knt,\mrt)$. The gradient is then time-dependent as well and given by \begin{align*}
    d_{\sn} \langle e(\tV,\tilde{\xi}),\psi \rangle  
    &= \int_{\tilde{\Omega}}  \tV \Bigl[ - \snt \nu \varphi_{\nu \nu} - 2 \snt \varphi_{\nu} -\rhot \varphi_x - \rhot \nu \varphi_{x \nu}\Bigr]\, dz\, d\tau\\
    d_{\rho} \langle e(\tV,\tilde{\xi}),\psi \rangle  
    &=   \int_{\tilde{\Omega}} \tV \Bigl[ - \snt \varphi_x - \snt \nu \varphi_{x \nu}\Bigr]\,dz\, d\tau, \\
    d_{\kn} \langle e(\tV,\tilde{\xi}),\psi \rangle 
    &=\int_{\tilde{\Omega}} \tV \Bigl[ (\mrt-\nu)\varphi_{\nu}-\varphi \Bigr] \,dz\, d\tau, \\
     d_{\mr} \langle e(\tV,\tilde{\xi}),\psi \rangle 
     &= \int_{\tilde{\Omega}}  \knt \tV \varphi_{\nu} \, dz\,.
\end{align*}

\subsection{Gradient descent algorithm for the parameter calibration}
Solving the first-order optimality condition all at once is difficult due to the forward-backward structure. 
Therefore, we propose a gradient descent algorithm in the following.

For a given initial parameter set $\xi_0$, we can solve the state equation for the Heston model with constant control variable $\xi$ \eqref{eq: Heston transformed} or time dependent parameter \eqref{eq: Heston transformed time}. 
With the state solution at hand, we can compute the corresponding adjoint equation \eqref{eq: Heston AD} or \eqref{eq:HestonAD_noA} %
backwards in time. 
Then we have all the information available to compute the gradient and update the parameter set using a gradient step. The procedure is sketched in algorithm \ref{alg:GDHeston}.

\bigskip
\begin{algorithm}[H]\label{alg:GDHeston}
	\SetAlgoLined
	\KwResult{calibrated parameters for Heston model}
	initialize parameters $u_0$\;
	\While{$\| \mathrm{gradient} \| > \epsilon$}{
		solve \eqref{eq: Heston transformed} or \eqref{eq: Heston transformed time}\;
		solve \eqref{eq: Heston AD} or \eqref{eq:HestonAD_noA}\; %
		compute the gradient\;
		line search for step size\;
		update the parameter set\; 
	}
	\caption{Gradient descent method for Heston parameter calibration}
\end{algorithm}

\bigskip 
Since the parameter domain for $\kappa_\nu$, $\mu_\nu$, $\sigma_\nu$ and $\rho$ is restricted, as well as the constraint that the Feller condition has to be fulfilled, we use the \textit{projected Armijo rule} \cite{Troeltzsch}.
In the projected Armijo rule, we choose the maximum $\sigma_k\in\{1,1/2,1/4,\ldots\}$ for which
\begin{equation*}
 f\bigl(\mathcal{P}(\xi_k - \sigma_k \nabla f(\xi_k))\bigr) - f(\xi_k) \le -\frac{\gamma}{\sigma_k}\|\mathcal{P}(\xi_k-\sigma_k \nabla f(\xi_k))-\xi_k\|_2^2.
\end{equation*}
Here $\gamma\in(0,1)$ is a numerical constant, which is problem-dependent and typically chosen as $\gamma=10^{-4}$. We will use this value for the numerical results later on.

\section{Discretization}\label{sec: Discretization} 
In this section we introduce a closure boundary condition at $\nu=0$ for the Heston and its adjoint and perform a domain truncation to discretize the problem. 
Following the Fichera theory and assuming that the Feller condition holds, no analytical boundary condition needs to be imposed at $\nu=0$, neither for the Heston model nor for its adjoint, since the PDEs have a pure outflow boundary at this point. Nevertheless, we need a closure condition for this boundary for the discretization. 
We suggest to follow Heston's approach, as discussed in \cite{Kutik.2015, clevenhaus2023ECMI} and use the reduced hyperbolic formulation of the Heston PDE and its adjoint. At $\Gamma_c$, we obtain 
\begin{equation}
    \tV_{\tau}=r\tV_x+\kn \mr \tV_{\nu}-r\tV\\
\end{equation}
for the $\log$-transformed normalized PDE and similar for the adjoint 
\begin{equation}\
    \varphi_{\tau} + (\sigma_\nu^2-\kappa_\nu\mr) \varphi_{\nu} + (q-r+\sigma_\nu \rho) \varphi_x  + (\kn -r) \varphi = -K(K\tV-\Vd).
\end{equation}
Now, we perform a domain truncation to obtain a rectangular grid, instead of a semi-unbounded domain and introduce the grid points. 
We consider uniform meshes in each direction and obtain for the spatial directions $x_i=x_{\min}+i \Delta_x$ for $i=0,\ldots,N_x$ 
with $\Delta_x=(x_{\max}-x_{\min})/N_x$ 
and $\nu_j=j \Delta_\nu$ for $j=0,\ldots,N_\nu$ with $\Delta_\nu=\nu_{\max}/N_\nu$, 
as well as $\tau_k=k\,\Delta_\tau$ for $k=0,\ldots,N_\tau$ with $\Delta_\tau=T/N_\tau$ for the temporal direction.

Since this is a proof-of-concept, simple and well-known spatial and temporal discretization methods are used to illustrate our approach.
For the time discretization, we use the well-known \textit{alternating direction implicit} (ADI) method, in more detail the Hundsdorfer-Verwer scheme \cite{Hundsdorfer.2002}. 
This scheme is of second order for any choice of $\theta$, where $\theta$ is a measure of classification similar to the $\theta$-method, and is able to handle mixed derivative terms.
To present the ADI method, we use a general second order PDE formulation
\begin{equation}
    u_\tau + a_{11} u_{\nu\nu} +2 a_{12} u_{x\nu} + a_{22} u_{xx} + b_1 u_{\nu} +b_2 u_x + c u=0.
\end{equation}
Since the $\log$-transformed normalised Heston PDE and its adjoint are second-order PDEs and we further use the same discretization.

In a first step, we split the operator of the PDE into three operators
\begin{equation}
    \F(\tau)=\F_0(\tau)+\F_1(\tau)+\F_2(\tau),
    \label{eq: Opterator Splitting ADI}
\end{equation}
with
\begin{align*}
    \F_0(\tau) &= 2 a_{12} D_{x \nu} \\
    \F_1(\tau) &= b_2 D_x + a_{22} D_{xx} - \half c I\\
    \F_2(\tau) &= b_1 D_{\nu} +  a_{11} D_{\nu \nu} -\half c I,
\end{align*}
where $D_x$ describes the discretization matrix of the first derivative w.r.t.\ $x$ and accordingly $D_{xx}$ of the second derivative w.r.t.\ $x$, $D_\nu$ and $D_{\nu\nu}$ of the first and second derivative w.r.t.\ $\nu$, $I$ denotes the identity matrix. 
The discretization matrices are derived using central finite differences. 
Let $u^{i,j}_k\approx u(x_i,\nu_j,\tau_k)$ and simplify $u_k\approx u(x,\nu,\tau_k)$.
In each time step, we have to solve the following system of equations
\begin{equation}
    \begin{cases}
    Y_0=u_{k}+\Dtau \F(\tau_{k}) u_{k}, \\
    Y_1=Y_{0}+\theta \Dtau \bigl(\F_1(\tau_{k+1})Y_1 -\F_1(\tau_{k}) u_{k} \bigr), \\
    Y_2=Y_{x}+\theta \Dtau \bigl(\F_2(\tau_{k+1},Y_2)-\F_2(\tau_{k}) u_{k} \bigr), \\
    \tY_0=Y_0+\half \Dtau \bigl(\F(\tau_{k+1})Y_2-\F(\tau_{k}) u_{k}\bigr),\\
    \tY_1=\tY_{0}+\theta \Dtau \bigl(\F_1(\tau_{k+1})\tY_1-\F_1(\tau_{k})u_{k}\bigr), \\
    \tY_2=\tY_{x}+\theta \Dtau \bigl(\F_2(\tau_{k+1})\tY_2-\F_2(\tau_{k})u_{k}\bigr), \\
    u_{k+1}=\tY_2.
    \end{cases}
\end{equation}
We choose $\theta=3/4$ and improve the implementation as we used the approach in \cite{Teng.2019} by using a matrix based instead of a vector based implementation of $u_k$.  

At this point, we have to discuss the boundary conditions for the Heston model again. We start with the $x$ dimension. 
We set $\tV(x_{\min},\nu,\tau)=\exp(-r\tau)$ and similar $\tV(x_{\max},\nu,\tau)=0$, as we suggest a sufficiently small $x_{\min}$ and a sufficiently large $x_{\max}$. 
For the variance boundaries, we follow the approach of K\`utik and Mikula \cite{Kutik.2015}.
Due to the Fichera theory, at $\nu=0$ we gain a outflow boundary and no information can enter the domain from the region $\nu<0$. 
Since the same holds for the adjoint, we propose the same approach there. In these, we extend the numerical domain for this case and introduce ghost cells $p_{i,0}=(x_i,\nu_0)$, $i=1,\ldots,N$ where $\nu_0 = -\Delta_y$ and determine the ghost cell values $u_{p_{i,0}}^n$, $i=1,\ldots,N$ by zero order extrapolation. 
We obtain a constant function
\begin{equation}
    u_{p_{i,0}}^n=u_{p_{i,1}}^n,\quad  i=1,\ldots,N.
\end{equation}
Now we consider the boundary, where $\nu \to \infty$ and the truncation shrinks the domain from $(0,\infty) \to (0,\nu_{\max})$. 
Dirichlet boundary condition from Heston would causes an unnatural jump in the solution. 
Therefore different strategies have been developed to overcome this issue, e.g.\ \cite{Kutik.2015}. 
Again, we follow K\`utik and Mikula and impose artificial homogeneous Neumann boundary conditions.
The choice is motivated by the variance independence of the original boundary condition from Heston, for sufficiently large $\nu_{\max}=\mathcal{O}(1)$.
Again, we use zero order extrapolation to implement this condition. 
We add ghost cells $p_{i,N_{\nu}+1}$, $i=1,\ldots,N$ where $\nu_{N_{\nu}+1} = \nu_{\max}+\Delta_y$ and determine their values $u_{p_{i,N_{\nu}+1}}^n$, $i=1,\ldots,N$ by a constant function
\begin{equation}
    u_{p_{i,N_{\nu}+1}}^n=u_{p_{i,N_{\nu}}}^n,\quad  i=1,\ldots,N.
\end{equation}
For a discussion about different numerical boundary conditions for the variance of the Heston model, we refer to \cite{clevenhaus2023ECMI}.
The integrals appearing in the gradient are computed by the trapezoidal rule.

\section{Numerical Results}\label{sec: Numerical Results}

Following \cite{clevenhaus2023ECMI}, we use 
\begin{equation}
    K=1.0, r=0.1, \xi_{\rm{ref}}=(5.0, 0.07, 0.5, -0.5)
\end{equation}
for generating an artificial $\Vd$ for each time step $\tau_k$. For the discretization, we use 
 $   N_{x}=79 , N_{\nu}=39, N_{\tau}=59.$
As bounds for the projected Armijo-rule for $\xi=(\kn,\mr,\sn,\rho)$, we set
 $   0 < \kn <8, \quad 0 < \mr <1, \quad 0< \sn < 1, \quad -1 < \rho < 1.$
Note that the projected Armijo-rule ensures that the Feller condition holds within each optimization step. Hence we are in the case of an outflow boundary. 
We set the maximal iteration value for the calibration to 20.
For the initial guesses $\xi_{\rm{init}}=(\kappa_{\nu}^{\rm{init}},\mu_{\nu}^{\rm{init}},\sigma_{\nu}^{\rm{init}},\rho^{\rm{init}})$ we used generated random numbers within a maximal percentage difference from $\xi_{\rm{ref}}$. 
We use four different percentages $10,25,50$ and $75$ and for each we generate five sets.
The initial parameters as well as the calibrated parameters in the constant case are given in Table~\ref{tab:init_kappa} for $\kappa_{\nu}$, Table~\ref{tab:init_mu} for $\mu_{\nu}$, Table~\ref{tab:init_sigma} for $\sigma_{\nu}$ and Table~\ref{tab:init_rho} for $\rho$. 

\begin{table}[]
    \centering
    \begin{tabular}{|c|c c c c c || c |c c c c c |}
    \hline
       $\kappa_{\nu}^{\rm{init}}$  & T1 & T2 & T3 & T4 & T5 & $\kappa_{\nu}^{\rm{cal}}$  & T1 & T2 & T3 & T4 & T5 \\ \hline
        10 & 4.60 & 5.15 & 4.50 & 4.70 & 4.50 &
        10 & 4.60 & 5.15 & 4.51 & 4.70 & 4.50
        \\
        25 & 4.75 & 4.60 & 5.25 & 5.90 & 4.45 &
        25 & 4.75 & 4.60 & 5.24 & 5.88 & 4.45
        \\
        50 & 4.10 & 3.25 & 3.55 & 2.90 & 7.00 &
        50 & 4.11 & 3.26 & 3.55 & 2.94 & 6.95
        \\
        75 & 5.70 & 7.30 & 7.00 & 7.50 & 6.60 &
        75 & 5.69 & 7.29 & 6.85 & 7.40 & 6.60
        \\ \hline
    \end{tabular}
    \caption{Five initial values sets for $\kappa_{\nu}^{\rm{init}}$ for the different percentages and the corresponding calibrated values for C0.}
    \label{tab:init_kappa}
\end{table}

\begin{table}[]
    \centering
    \footnotesize
    \begin{tabular}{|c|c c c c c|| c| c c c c c|}
    \hline
        $\mu_{\nu}^{\rm{init}}$  & T1 & T2 & T3 & T4 & T5 &
        $\mu_{\nu}^{\rm{cal}}$  & T1 & T2 & T3 & T4 & T5 
        \\ \hline
        10 & 0.0679 & 0.0637 & 0.0735 & 0.0714 & 0.0714&
        10 & 0.0598 & 0.0727 & 0.0575 & 0.0632 & 0.0577
        \\
        25 & 0.0721 & 0.0777 & 0.0616 & 0.0630 & 0.0868 &
        25 & 0.0655 & 0.0610 & 0.0720 & 0.0836 & 0.0615
        \\
        50 & 0.1001 & 0.0406 & 0.0840 & 0.0448 & 0.0434 &
        50 & 0.0476 & 0.0148 & 0.0484 & 0.0136 & 0.0882
        \\
        75 & 0.0742 & 0.0532 & 0.0448 & 0.0399 & 0.0490 &
        75 & 0.0847 & 0.1069 & 0.1005 & 0.1064 & 0.0928
        \\ \hline
    \end{tabular}
    \caption{Five initial values sets for $\mu_{\nu}^{\rm{init}}$ for the different percentages and the corresponding calibrated values for C0.}
    \label{tab:init_mu}
\end{table}

\begin{table}[]
    \centering
    \begin{tabular}{|c|c c c c c|| c| c c c c c|}
    \hline
        $\sigma_{\nu}^{\rm{init}}$  & T1 & T2 & T3 & T4 & T5 &
        $\sigma_{\nu}^{\rm{cal}}$  & T1 & T2 & T3 & T4 & T5 
        \\ \hline
        10 & 0.465& 0.545& 0.510& 0.515& 0.480&
        10 & 0.451& 0.559& 0.484& 0.502& 0.455
        \\
        25 & 0.575& 0.570& 0.505& 0.415& 0.480&
        25 & 0.565& 0.546& 0.523& 0.455& 0.429
        \\
        50 & 0.580& 0.360& 0.630& 0.345& 0.655&
        50 & 0.501& 0.292& 0.581& 0.261& 0.713
        \\
        75 & 0.705& 0.470& 0.235& 0.350& 0.565&
        75 & 0.718& 0.562& 0.390& 0.489& 0.630
        \\ \hline
    \end{tabular}
    \caption{Five initial values sets for $\sigma_{\nu}^{\rm{init}}$ for the different percentages and the corresponding calibrated values for C0.}
    \label{tab:init_sigma}
\end{table}

\begin{table}[]
    \centering
    \footnotesize
    \begin{tabular}{|c|c c c c c|| c| c c c c c|}
    \hline
       $\rho^{\rm{init}}$  & T1 & T2 & T3 & T4 & T5&
       $\rho^{\rm{cal}}$  & T1 & T2 & T3 & T4 & T5\\ \hline
        10 & -0.550 & -0.485 & -0.520 & -0.520 & -0.495 &
        10 & -0.483 & -0.547 & -0.404 & -0.460 &  -0.384
        \\
        25 & -0.410 & -0.490 & -0.605 & -0.585 & -0.380 &
        25 & -0.370 & -0.390 & -0.690 & -0.829 & -0.159
        \\
        50 & -0.555 & -0.465 & -0.350 & -0.685 & -0.725 &
        50 & -0.258 & -0.242 & -0.207 & -0.441 & -0.847
        \\
        75 & -0.695 & -0.205 & -0.775 & -0.655 & -0.150 &
        75 & -0.743 & -0.972 & -0.990 & -0.990 & -0.437
        \\ \hline
    \end{tabular}
    \caption{Five initial values sets for $\rho^{\rm{init}}$ for the different percentages and the corresponding calibrated values for C0.}
    \label{tab:init_rho}
\end{table}

We observe that the calibration leaves $\kn$ nearly untouched, whereas the other parameter values are significantly changed. %
This could be reasoned by the structure of the drift term $\kn(\mr-\nu)$, since $\kn$ and $\mr$ are multiplied.
As a optimization measure, we compute the relative reduction of the cost functional using
\begin{equation} \label{eq: rel_red}
    r(\xi_{\rm{init}})= 100 \cdot \left(1-\frac{J(\tV,\xi_{\rm opt})}{J(\tV,\xi_{\rm{init}})}\right). 
\end{equation}
Table \ref{tab:rel_C0} shows the relative reduction of the cost functional for the different test cases for $\xi_{\rm{init}}$.

\begin{table}[]
    \centering
    \begin{tabular}{|c|c c c c c|} \hline
         C0  & T1 & T2 & T3 & T4 & T5 \\ \hline
        10 &43.12 & 77.38 & 65.22 & 63.22 & 52.75\\
        25 &58.46 & 78.38 & 62.14 & 66.06 & 64.53\\
        50 &83.09 & 20.84 & 52.76 & 30.98 & 72.62\\
        75 &16.23 & 72.89 & 82.89 & 82.09 & 80.43\\ \hline
        \end{tabular}
    \caption{Relative reduction of the cost functional computed with \ref{eq: rel_red} using the different test cases for $\xi_{\rm{init}}$ and $\xi_{\rm{cal}}$ in the constant calibration setting.}
    \label{tab:rel_C0}
\end{table}
\begin{figure}
    \centering
    \includegraphics[scale=.4]{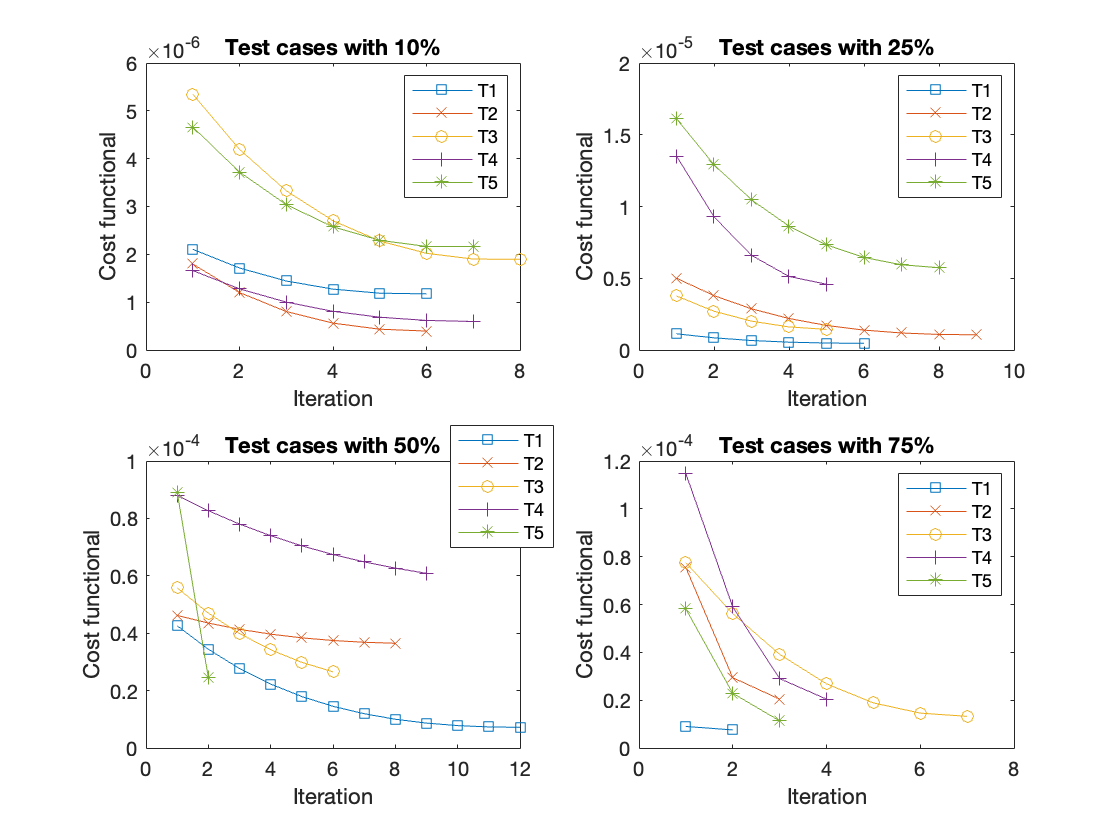}
    \caption{Cost functional evolution per iteration for the test cases within the constant parameter calibration.}
    \label{fig:C0}
\end{figure}

As the calibrated values differ over the test cases, it is reasonable to assume, that we only find local minima. 
However the results are remarkable.%

Since in the real market the parameter are not considered as constant, we improve the approach by considering different parameters, and some parameter sets as time-dependent. 
From the relative change in the constant calibration setting, we choose the following (additional) test cases listed in Table~\ref{tab:cases}. 
Further the table includes the links to the cost functional reduction tables and figures of the cost function evolution with respect to the different cases.

\begin{table}[]
    \centering
    \small
    \begin{tabular}{|c|c| c| c| c|}
    \hline
           &  $\kn$ & $\mr$ & $\sn$ & $\rho$  \\ \hline
       C0  & constant & constant& constant& constant  \\ 
       C1  &  time-dependent & constant & constant & constant\\
       C2  &  constant & time-dependent & constant & constant \\
       C3  &  constant & constant & time-dependent & constant \\
       C4  &  constant & constant & constant & time-dependent \\
       C5  &  time-dependent & time-dependent & time-dependent & time-dependent \\
       C6  &  constant & time-dependent & time-dependent & time-dependent \\
       C7  &  constant & time-dependent & constant & time-dependent \\ \hline
    \end{tabular}
    \caption{Different cases for calibration setting.}
    \label{tab:cases}
\end{table}
For each case of the time-dependent test cases $\Vd$ is generated as before and $\xi^{\rm{init}}$ is assumed to be constant and thus also the same as before. 
\begin{figure}
    \centering
    \includegraphics[scale=.4]{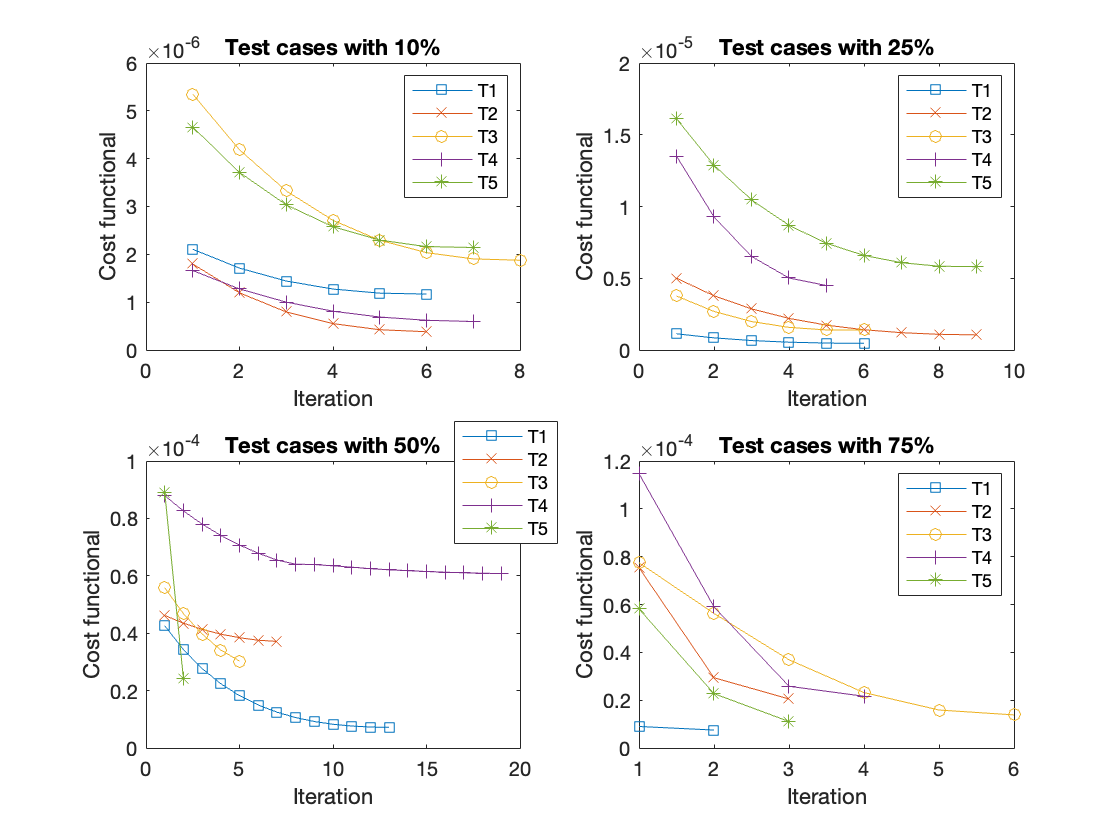}
    \caption{Cost functional evolution per iteration for the different test set in the case scenario C1.}
    \label{fig:C1}
\end{figure}
\begin{figure}
    \centering
    \includegraphics[scale=.4]{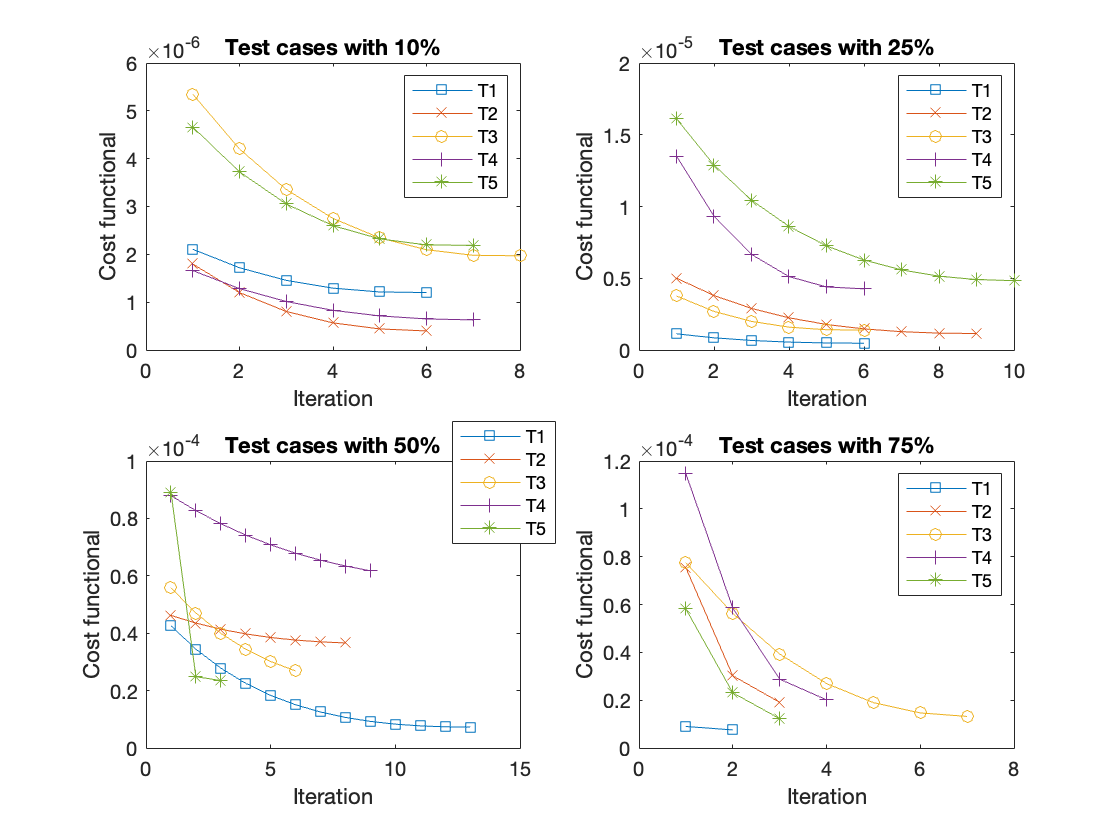}
    \caption{Cost functional evolution per iteration for the different test set in the case scenario C2.}
    \label{fig:C2}
\end{figure}
\begin{figure}
    \centering
    \includegraphics[scale=.4]{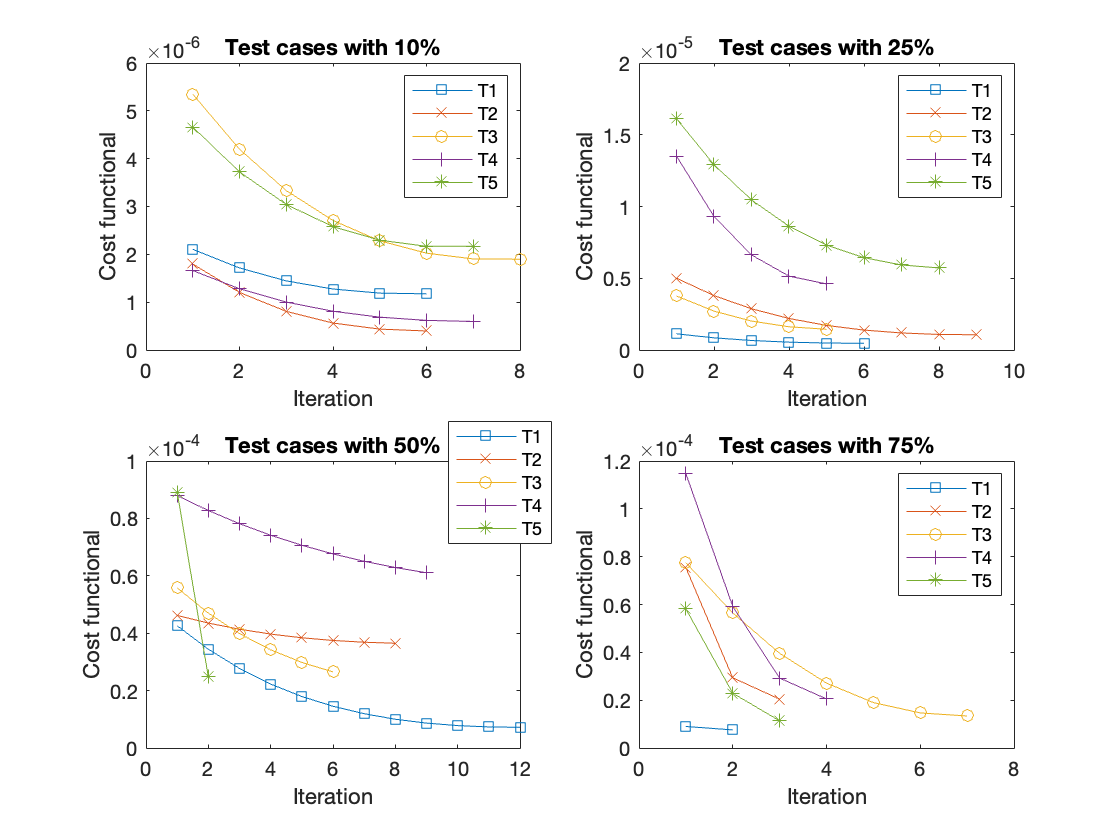}
    \caption{Cost functional evolution per iteration for the different test set in the case scenario C3.}
    \label{fig:C3}
\end{figure}
\begin{figure}
    \centering
    \includegraphics[scale=.4]{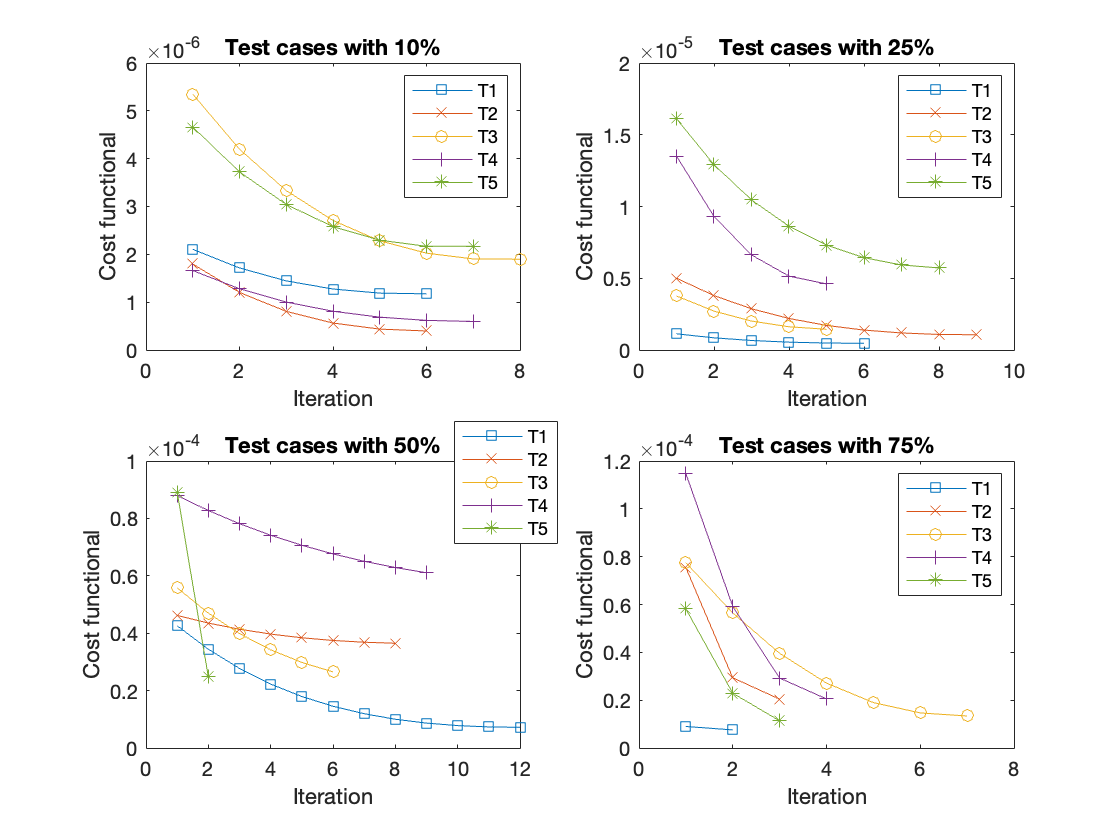}
    \caption{Cost functional evolution per iteration for the different test set in the case scenario C4.}
    \label{fig:C4}
\end{figure}
\begin{figure}
    \centering
    \includegraphics[scale=.4]{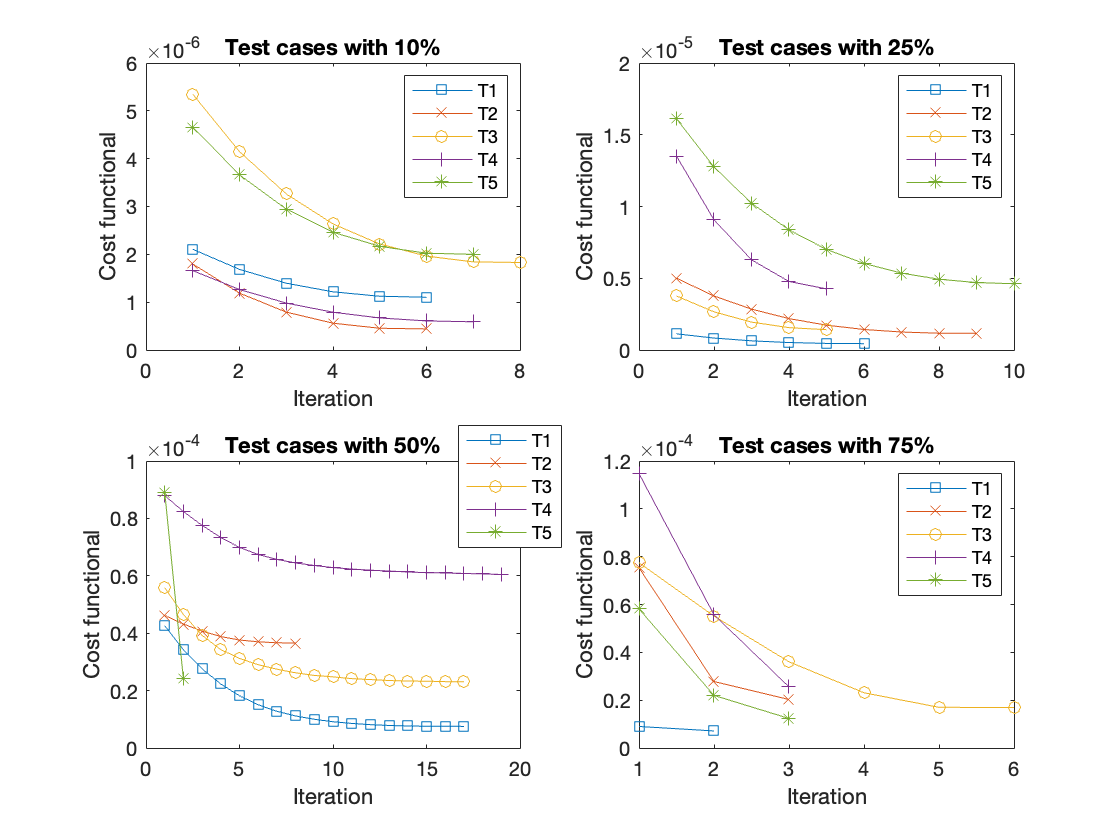}
    \caption{Cost functional evolution per iteration for the different test set in the case scenario C5.}
    \label{fig:C5}
\end{figure}
\begin{figure}
    \centering
    \includegraphics[scale=.4]{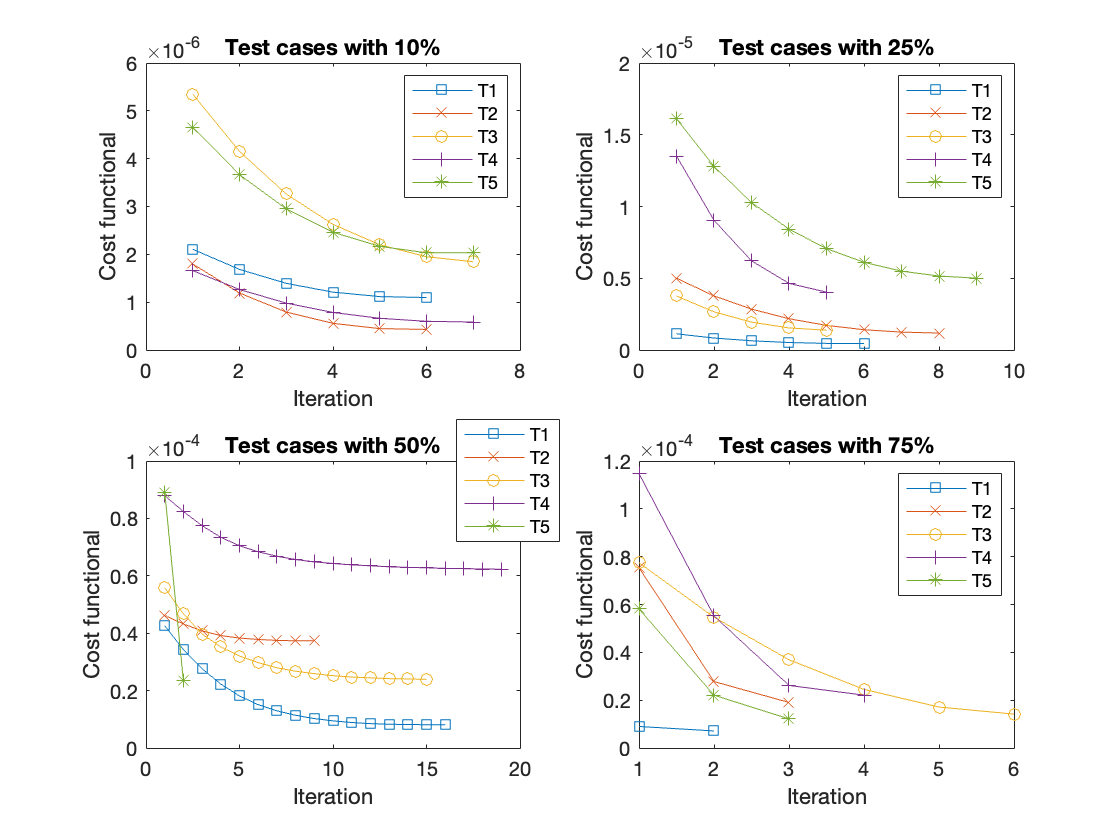}
    \caption{Cost functional evolution per iteration for the different test set in the case scenario C6.}
    \label{fig:C6}
\end{figure}
\begin{figure}
    \centering
    \includegraphics[scale=.4]{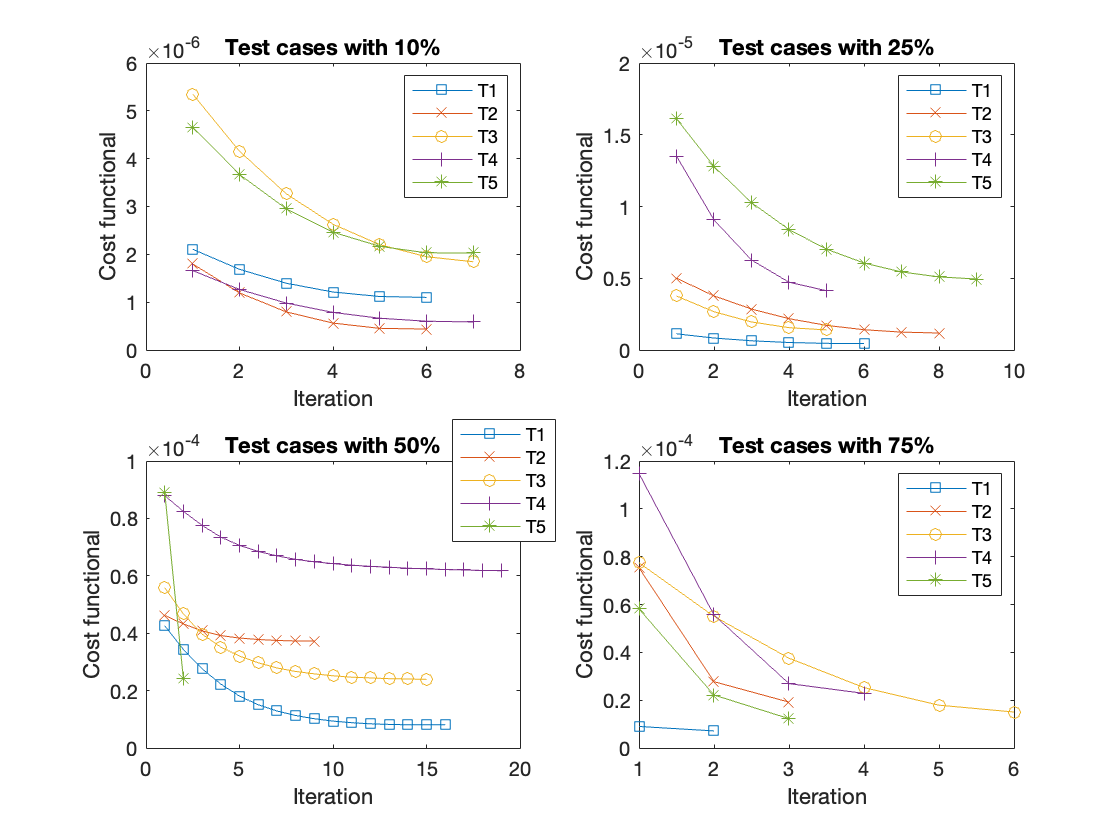}
    \caption{Cost functional evolution per iteration for the different test set in the case scenario C7.}
    \label{fig:C7}
\end{figure}

\begin{table}[]
    \centering
    \begin{tabular}{|c|c c c c c|} \hline 
        C1&		T1 & T2 & T3 & T4 & T5 \\	\hline				
        10 &    44.48&	78.75&	64.91&	64.12&	54.03\\
        25 &	58.97&	78.98&	62.96&	66.77&	63.96\\
        50 &	83.13&	19.59&	45.71&	33.45&	72.84\\
        75 &	16.59&	72.80&	82.14&	81.13&	80.94\\ \hline 
         \end{tabular}
    \caption{Relative cost function reduction for the C1 calibration.}
    \label{tab:rel_C1}
\end{table}
\begin{table}[]
    \centering
    \begin{tabular}{|c|c c c c c|} \hline 
        C2&		T1 & T2 & T3 & T4 & T5 \\		\hline			
        10 &    42.96&	77.85&	63.17&	62.09&	53.05\\
        25 &	58.03&	77.05&	63.09&	68.46&	70.00\\
        50 &	82.90&	20.79&	51.89&	29.79&	73.66\\
        75 &	16.05&	74.64&	82.95&	82.43&	78.89\\ \hline 
         \end{tabular}
    \caption{Relative cost function reduction for the C2 calibration.}
    \label{tab:rel_C2}
\end{table}
\begin{table}[]
    \centering
    \begin{tabular}{|c|c c c c c|} \hline 
        C3&		T1 & T2 & T3 & T4 & T5 \\		\hline			
        10 &    44.22&	77.89&	64.50&	64.04&	53.52\\
        25 &	58.49&	78.76&	61.62&	65.92&	64.50\\
        50 &	83.05&	20.97&	52.49&	30.55&	72.23\\
        75 &	15.96&	72.91&	82.71&	82.10&	80.44\\ \hline 
         \end{tabular}
            \caption{Relative cost function reduction for the C3 calibration.}
    \label{tab:rel_C3}
\end{table}
\begin{table}[]
    \centering
    \begin{tabular}{|c|c c c c c|} \hline 
        C4&		T1 & T2 & T3 & T4 & T5 \\			\hline		
        10 &    49.03&	75.90&	66.62&	66.49&	56.81\\
        25 &	60.42&	77.71&	62.64&	68.50&	62.96\\
        50 &	78.54&	18.24&	52.63&	30.29&	74.03\\
        75 &	20.89&	73.59&	81.75&	80.49&	79.74\\ \hline 
    \end{tabular}
    \caption{Relative cost function reduction for the C4 calibration.}
    \label{tab:rel_C4}
\end{table}
\begin{table}[]
    \centering
    \begin{tabular}{|c|c c c c c|} \hline 
        C5&		T1 & T2 & T3 & T4 & T5 \\		\hline			
        10 &    47.63&	75.50&	65.79&	64.45&	57.08\\
        25 &	60.70&	76.67&	62.36&	68.51&	71.33\\
        50 &	82.32&	20.89&	58.63&	31.72&	72.96\\
        75 &	20.79&	73.10&	78.14&	77.71&	78.76\\ \hline 
    \end{tabular}
    \caption{Relative cost function reduction for the C5 calibration.}
    \label{tab:rel_C5}
\end{table}
\begin{table}[]
    \centering
    \begin{tabular}{|c|c c c c c|} \hline 
        C6&		T1 & T2 & T3 & T4 & T5 \\			\hline		
        10 &    47.91&	76.08&	65.50&	64.77&	56.43\\
        25 &	60.27&	76.39&	63.37&	70.16&	68.91\\
        50 &	80.94&	19.19&	57.15&	29.38&	73.71\\
        75 &	20.84&	74.74&	81.70&	80.74&	79.10\\ \hline 
    \end{tabular}
    \caption{Relative cost function reduction for the C6 calibration.}
    \label{tab:rel_C6}
\end{table}
\begin{table}[]
    \centering
    \begin{tabular}{|c|c c c c c|} \hline 
        C7&		T1 & T2 & T3 & T4 & T5 \\			\hline		
        10 &    47.79&	75.66&	65.50&	64.66&	56.48\\
        25 &	60.21&	76.37&	62.80&	69.43&	69.29\\
        50 &	81.08&	19.35&	57.16&	30.05&	73.02\\
        75 &	20.62&	74.55&	80.59&	80.09&	78.89\\ \hline
    \end{tabular}
    \caption{Relative cost function reduction for the C7 calibration.}
    \label{tab:rel_C7}
\end{table}

To quantify the different calibrations, we compute the average relative cost functional reduction by
\begin{equation}
    r_a(\xi_{\rm{init}})=\frac{r(\xi_{\rm{init}})}{20}
\end{equation}
and summarize the results in Table~\ref{tab:res_all}. Note, that all cost function reduction averages are huge. 
We observe that a time-dependent calibration for one parameter alone (C1-C4) doesn't improve the cost functional reduction significantly. The first slightly improvement can be found by using at least two time-dependent parameters (C5-C7). 
Surprisingly C7 gives the best calibration results, where $\kn$ is the only variable which is calibrated as a constant, and C5, where all parameters are calibrated as time-dependent gives the least improvement, when considering combinations of time-dependent parameter calibration.
The fact that $\Vd$ is generated with constant parameters and the best case considers time-dependent parameters indicates that the time-dependency is a good way to overcome the local minima. 
Those results are inline with literature \cite{Mikhailov.2003,Teng.2016}.
\begin{table}[]
    \centering
    \small
    \begin{tabular}{|c|c|c|c|c|c|c|c|c|} \hline
        Case Setup & C0 & C1 & C2& C3& C4& C5& C6& C7\\ \hline
       $r_a(\xi_{\rm{init}})$ & 61.30 & 61.31 & 61.49 & 61.34 & 61.86 & 62.25 & 62.36 & 66.12 \\ 
       table & \ref{tab:rel_C0} & \ref{tab:rel_C1} & \ref{tab:rel_C2} & \ref{tab:rel_C3} & \ref{tab:rel_C4} & \ref{tab:rel_C5} & \ref{tab:rel_C6} & \ref{tab:rel_C7}\\
       figure & \ref{fig:C0} & \ref{fig:C1} & \ref{fig:C2} & \ref{fig:C3} & \ref{fig:C4} & \ref{fig:C5} & \ref{fig:C6} & \ref{fig:C7}\\ \hline
    \end{tabular}
    \caption{Average cost functional reduction in percentage for the different cases over the 20 test cases, as well as a list for the corresponding relative cost functional reduction as well as the figures of the cost functional evolution.}
    \label{tab:res_all}
\end{table}

\section{Conclusion and Outlook}\label{sec: Conclusion and Outlook}
Our paper began with an introduction to the Heston model. 
The introduction was followed by the derivation of a gradient-based optimization, including the derivation of a gradient descent algorithm. 
In the next section, the discretization of the schemes and algorithms is presented and numerical results are illustrated.
These show that the gradient descent algorithm is a feasible calibration technique. 
The relative cost function reduction is remarkable, even though we can expect to find only local minima. 
Thus, our approach follows recent research. 
Furthermore, the assumption of at least time-dependent parameters is better suited to the real market situation, since in the real market almost no parameter is constant. 
In future work, we plan to incorporate the gradient descent algorithm into the space mapping approach and to test the space mapping approach on real market data.

\bibliographystyle{plain}
\bibliography{Paper22}

\end{document}